\documentclass[journal]{IEEEtran}

\usepackage{amsmath,amssymb}
\usepackage{amsfonts}
\usepackage{mathrsfs}
\usepackage{color}
\usepackage{graphicx}
\usepackage{epsfig} 

\usepackage{array,booktabs} 
\usepackage{amsmath,amssymb} 
\usepackage{stfloats} 
\usepackage{psfrag} 
\usepackage{color}

\usepackage{enumitem}

\definecolor{box_color}{rgb}{.8,.8,.8}

\usepackage[utf8]{inputenc}
\usepackage[T1]{fontenc}
\usepackage{lmodern}
\usepackage{epstopdf}
\newtheorem{lemma}{Lemma}

\newtheorem{corollary}{Corollary}
\newtheorem{fact}{Fact}

\newtheorem{remark}{Remark}

\newtheorem{assumption}{Assumption}
\newtheorem{ass}{C.}

\def\begcen{\begin{center}}
\def\endcen{\end{center}}

\newcommand{\RE}{\mathbb {R}}    

\def\L2{{\cal L}_2}
\def\L2e{{\cal L}_{2e}}

\def\rea{\mathbb{R}}

\def\begequarr{\begin{eqnarray}}
\def\endequarr{\end{eqnarray}}
\def\begequarrs{\begin{eqnarray*}}
\def\endequarrs{\end{eqnarray*}}
\def\begarr{\begin{array}}
\def\endarr{\end{array}}
\def\begequ{\begin{equation}}
\def\endequ{\end{equation}}
\def\lab{\label}
\def\begdes{\begin{description}}
\def\enddes{\end{description}}
\def\begenu{\begin{enumerate}}
\def\begite{\begin{itemize}}
\def\endite{\end{itemize}}
\def\begmyite{\begin{myitemize}}
\def\endmyite{\end{myitemize}}
\def\endenu{\end{enumerate}}

\def\lef[{\left[\begin{array}}
\def\rig]{\end{array}\right]}

\def\begcen{\begin{center}}
\def\endcen{\end{center}}
\def\begrem{\begin{remark}\rm}
\def\endrem{\end{remark}}
\def\begassum{\begin{assumption}}
\def\endassum{\end{assumption}}
\def\begassums{\begin{assumption*}}
\def\endassums{\end{assumption*}}
\def\begassu{\begin{ass}}
\def\endassu{\end{ass}}
\def\beglem{\begin{lemma}}
\def\endlem{\end{lemma}}
\def\begcor{\begin{corollary}}
\def\endcor{\end{corollary}}
\def\begfac{\begin{fact}}
\def\endfac{\end{fact}}

\def\begmat#1{\begin{bmatrix}#1\end{bmatrix}}
\def\begali#1{\begin{align}{#1}\end{align}}
\def\begalis#1{\begin{align*}{#1}\end{align*}}

\def\calf{{\cal F}}

\def\calf{{\cal F}}

\def\L2e{{\cal L}_{2e}}

\def\rea{\mathbb{R}}

\usepackage{color}

\usepackage[prependcaption,colorinlistoftodos]{todonotes}

\begin{document}

\title{ Correction to the paper ``A robust IDA-PBC approach for handling uncertainties in underactuated mechanical systems''  }
\author{Alejandro~Donaire,~\IEEEmembership{Member,~IEEE,}
        Jose~Guadalupe~Romero,
        Romeo~Ortega,~\IEEEmembership{Fellow,~IEEE,}
\thanks{A. Donaire is with the School of Engineering, The University of Newcastle, University Dr, Callaghan, 2308, NSW, Australia, e-mail: alejandro.donaire@newcastle.edu.au.}
\thanks{J. G. Romero is with the Departamento Acad\'{e}mico de Sistemas Digitales, ITAM, R\'{i}o Hondo No.1, 01080, Ciudad de M\'exico, M\'{e}xico, e-mail:jose.romerovelazquez@itam.mx}
\thanks{R. Ortega is with Laboratoire des Signaux et Syst\`emes, CNRS--SUPELEC, Plateau du Moulon, 91192
Gif--sur--Yvette, France, e-mail: ortega@lss.supelec.fr}
}

\maketitle

\begin{abstract}
In this note, it is shown that the results claimed in the paper \cite{Ryalat2018}---as well as the examples presented there---are, unfortunately, incorrect. 
\end{abstract}
\IEEEpeerreviewmaketitle

\noindent {\bf Notation.} $I_n$ is the $n \times n$ identity matrix and $0_{n \times s}$ is an
$n \times s$ matrix of zeros, $0_n$ is an $n$--dimensional column vector of zeros. For $x \in \rea^n$, $S \in \rea^{n \times n}$, $S=S^\top>0$, we denote the Euclidean norm $|x|^2:=x^\top x$, and the weighted--norm $\|x\|^2_S:=x^\top S x$.  We use the notation \cite[\hspace{-1mm}$(\#)$]{Ryalat2018} to refer to the equation number $(\#)$ in  \cite{Ryalat2018}. To simplify the notation, the arguments of the various mappings are indicated only the first time they are defined. 
%
\section{Background}
\lab{sec1}
To set-up the notation we briefly review in this section the interconnection and damping assignment passivity-based control (IDA-PBC) method proposed in \cite{ORTtac} for underactuated mechanical systems, which are described in port-Hamiltonian (pH) form by 
\begin{eqnarray}
\left[ \begin{array}{c}
                   \dot q \\
                   \dot p
                   \end{array} \right] = \left[ \begin{array}{cc}
                                               0_{n\times n}      & I_{n} \\
                                              -I_{n}  & 0_{n\times n}
                                              \end{array} \right] \nabla H(q,p)+ \left[ \begin{array}{c}
                                                                        0_{n\times m} \\
                                                                        G(q)
                                                                  \end{array}\right] \, u,
\label{sys}
\end{eqnarray}
\noindent where $q,p \in \rea^{n}$ are the generalized position and momenta, respectively, $u \in \rea^{m},\; m < n$ is the control, $G \colon \rea^{n} \to \rea^{n \times m}$, with $\mbox{\rm rank}(G)= m$. The function $H \colon \rea^{n} \times \rea^{n} \to \rea,$
\begin{equation*}
H(q,p) := {1 \over 2} \, p^\top \,  M^{-1}(q) \, p + V(q)
\end{equation*}
is the total energy with $M \colon \rea^{n}  \to \rea^{n\times n}$, the positive definite mass matrix and $V \colon \rea^{n} \to \rea$ the potential energy.

The control objective is to design a static, state-feedback controller that assigns to the closed loop a desired stable equilibrium  $(q,p) = (q^{\star},0_n)$, $q^{\star} \in \rea^{n}$. This is achieved in IDA-PBC by assigning to the closed loop the pH target dynamics
\begali{
\lef[{c}
             \dot{q} \\
             \dot{p}
             \rig] & = \calf_d(q,p)  \nabla H_d(q,p),
\label{sysd} \\ \nonumber
\calf_d(q,p) \hspace{-.5mm}&:=\hspace{-.5mm} \lef[{cc} 
                      0_{n\times n}            &  M^{-1}(q) \, M_d(q)   \\
                     \hspace{-1.5mm}-M_d(q) \, M^{-1}(q)  &  \hspace{-1.5mm}J_2(q,p)\hspace{-.5mm}-\hspace{-.5mm}G(q)K_P G^\top(q)
                     \rig], 
}                    
with new total energy function $H_{d} \colon \rea^{n} \times \rea^{n} \to \rea,$
\begin{equation}
\lab{hd}
H_{d}(q,p) := {1 \over 2} \, p^\top \, M_{d}^{-1}(q)\, p + V_{d}(q),
\end{equation}
where the desired mass matrix $M_{d} \colon \rea^{n}  \to \rea^{n \times n}$ is positive definite,  the desired potential energy {\bf $V_{d} \colon \rea^{n} \to \rea$} verifies
\begin{equation*}
q^{\star}   =  \arg \min V_{d}(q),
\end{equation*}
and  $K_P \in \rea^{m \times m}$ is a free positive definite matrix. The matrix $J_{2} \colon \rea^{n} \times \rea^{n} \to \rea^{n \times n}$  is free to the designer and fulfills the skew-symmetry condition
\begin{equation*}
J_{2}(q,p)  =  -J_{2}^\top(q,p).
\end{equation*}
The pH target dynamics of the closed loop \eqref{sysd} ensures the following properties \cite{ORTtac}: 
\begenu
\item[{\bf S1}]  The closed-loop system  (\ref{sysd}) has a stable equilibrium point at $(q^{\star},0)$ with Lyapunov function $H_{d}$, which verifies 
$$
\dot{H}_{d} = -\| G^\top M_d^{-1}p\|_{K_P}^2 \leq 0.
$$
\item[{\bf S2}]  The equilibrium is asymptotically stable provided that the output 
\begin{equation*}
y_d:= G^\top M_d^{-1}p
\end{equation*}
is detectable {\em with respect to the dynamics} \eqref{sysd}.
\endenu

By equating the right-hand sides of \eqref{sys} and \eqref{sysd} one obtains the so-called {\em matching equations}, which are two partial differential equations (PDEs) that identify the assignable $M_d$ and $V_d$, and gives an explicit expression for the  (static state-feedback) control signal $u$.

\section{Formulation of the Robust IDA-PBC Problem in \cite{Ryalat2018}}
\lab{sec2}
%
The authors of \cite{Ryalat2018} consider that the system \eqref{sys} is in closed loop with an IDA-PBC yielding the closed-loop dynamics \eqref{sysd}. Moreover, it is assumed---{\em i.e.}, Assumption 3.1 of \cite{Ryalat2018}---that the equilibrium  $(q,p)=(q^{\star},0_n)$ of  \eqref{sysd} is {\em asymptotically stable}. The control objective in \cite{Ryalat2018} is to robustify the IDA-PBC with respect to the presence of additive disturbances. 

We summarise here the problem formulation in \cite{Ryalat2018}.\\

\noindent \emph{Problem formulation.} Given the dynamics 
\begequ
\label{sysdist}
\left[ \begin{array}{c}  \dot{q} \\  \dot{p} \end{array} \right] = \calf_d\nabla H_d +  \left[ \begin{array}{c}   0_{n\times m} \\ G \end{array}\right] \, v  + \left[ \begin{array}{c} d_1 \\  d_2 \end{array}\right],
\endequ
verifying Assumption 3.1 of \cite{Ryalat2018}, where $d_1 \in \RE^n$ and $d_2 \in \RE^n$ are said to be the ``unmatched'' and ``matched'' disturbances, respectively. Find a  control law of the form
\begali{
\nonumber
\dot x_v &= F(q,p,x_v)\\
\nonumber
v&=\beta(q,p,x_v),
}
where $x_v \in \RE^n$, that ensures the following objectives:
\begenu
\item[{\bf O1.}] asymptotic stability of the equilibrium $(q,p,x_v)=(q^{\star},0_n,0_n)$ when $d_1=d_2=0$; 
\item[{\bf O2.}]  integral-input-to-state stability (iISS) or input-to-state stability (ISS) with respect to the disturbances of the closed-loop.
\endenu

\begrem
\lab{rem1} 
The controller proposed in \cite[\hspace{-1mm}(15)]{Ryalat2018} is called an {\em integral} control even though the right hand side of the dynamic extension $\dot x_v$ is clearly a function of $x_v$.
\endrem

\begrem
\lab{rem2} 
The qualification of $d_1$ and $d_2$ as unmatched and matched disturbances is arbitrary and has {\em no relation} with the null space or image of the input matrix $G$.
\endrem

\section{Main Claims in \cite{Ryalat2018} and their Rebuttal}
\lab{sec3}
The main claims of \cite{Ryalat2018} are contained in {Proposition 4.1}, {Section V.A} and {Proposition 5.1} of \cite{Ryalat2018}. Unfortunately, as shown below, the proofs are not correct and thus the claims are invalid.
\subsection{Proposition 4.1 of \cite{Ryalat2018}}
The claims in this proposition, which refer to the dynamics of the system \eqref{sysdist}, with $d_1=0$ and $d_2=0$, in closed-loop with the control law \cite[\hspace{-1mm}(15)]{Ryalat2018}, are the following.\\

\noindent {\bf C1.}  The closed-loop system can be written in the pH form \cite[\hspace{-1mm}(20)]{Ryalat2018}, which is described in the new coordinates defined in \cite[\hspace{-1mm}(17)]{Ryalat2018}. 

\noindent {\bf C2.} The control objective {\bf O1} is verified.
\\ \ \\
Both claims are, unfortunately, incorrect for the following reasons.\\

\noindent {\bf R1.} The dynamics of the closed-loop using control law \cite[\hspace{-1mm}(15)]{Ryalat2018} does not have the form \cite[\hspace{-1mm}(20)]{Ryalat2018} without further assumptions on the mechanical system.
The matching of the open and closed-loop dynamics results in \cite[\hspace{-1mm}(22)]{Ryalat2018}, which has the form
\begin{eqnarray}
G\,v &=&(J_2-R_d)M_d^{-1}\mathcal{K} x_v - \mathcal{K} \mathcal{K}^\top M^{-1} \nabla_qV_d+ \nonumber \\ 
&&+\frac 12 M_dM^{-1}\nabla_q(p^\top M_d^{-1}p) -\frac 12 M_dM^{-1}\nabla_q \nonumber\\
&&\times \left[ (p+\mathcal{K} x_v)^\top M_d^{-1} (p+\mathcal{K} x_v) \right],\nonumber
\end{eqnarray}
with $\mathcal K(q):= G(q)K_iG^\top(q)$, and $K_i =K_i^\top>0$. The control law, proposed in the first equation of \cite[\hspace{-1mm}(15)]{Ryalat2018}, is given by 
\begin{eqnarray}
v &=&(G^\top G)^{-1} G^\top \Big\{ (J_2-R_d)M_d^{-1}\mathcal{K} x_v - \mathcal{K} \mathcal{K}^\top M^{-1} \nabla_qV_d \nonumber \\
&& + \frac 12 M_dM^{-1} \nabla_q(p^\top M_d^{-1}p)  \nonumber\\
&& -\frac 12 M_dM^{-1}\nabla_q\left[ (p+\mathcal{K} x_v)^\top M_d^{-1} (p+\mathcal{K} x_v) \right]
\Big \}. \nonumber
\end{eqnarray}
However, to ensure that the matching is satisfied, one should verify that
\begin{eqnarray} \label{matching1}
&& G^\perp \Big\{ (J_2-R_d)M_d^{-1}\mathcal{K} x_v - \mathcal{K} \mathcal{K}^\top M^{-1} \nabla_qV_d   \nonumber \\
&&+\frac 12 M_dM^{-1} \nabla_q(p^\top M_d^{-1}p) \lab{gper} \\
\nonumber
&& -\frac 12 M_dM^{-1}\nabla_q\left[ (p+\mathcal{K} x_v)^\top M_d^{-1} (p+\mathcal{K} x_v) \right] 
\Big\} =0
\end{eqnarray}
is satisfied. The condition \eqref{matching1} is not verified for general underactuated mechanical systems and the claim {\bf C1} is false. 

\noindent {\bf R2.} Asymptotic stability is not ensured from Assumption 3.1. 

As indicated in Section \ref{sec2}, Assumption 3.1 of \cite{Ryalat2018} is that the equilibrium of \eqref{sysd} is asymptotically stable. It is argued in \cite{Ryalat2018} that the control objective \textbf{O1} ``is concluded using the arguments used in Assumption 3.1'', but it is not clear at all to which arguments the authors are referring to. In a rambling text in the paragraph where Assumption 3.1 is enunciated and the one below, reference is made to Proposition 1 in \cite{ORTtac}. This proposition simply quotes the statement {\bf S2} of Section \ref{sec1}, where it is important to underscore that the detectability property refers to the system \eqref{sysd}. Even assuming that the bound \cite[\hspace{-1mm}(21)]{Ryalat2018} is correct, which it is not because---as proven in {\bf R1}---the closed-loop dynamics is {\em not} given by \cite[\hspace{-1mm}(20)]{Ryalat2018}---it is erroneous to assume that the aforementioned detectability property with respect to \eqref{sysd} implies detectability with respect to \cite[\hspace{-1mm}(20)]{Ryalat2018}. The discussion that is given in the remaining lines of the paragraph below \cite[\hspace{-1mm}(21)]{Ryalat2018} have no connection with the asymptotic stability claim {\bf C2}.      \\

\begrem
\lab{rem3} 
The mistake made in Claim {\bf C1} stems from the following elementary linear algebra fact. The equation $Az=b$, with  $A \in \rea^{n \times m}$ a full-rank, tall matrix and  $b \in \rea^n$ admits the solution $z=(A^\top A)^{-1}A^\top b$ {\em if and only if} $A^\perp b=0$, with  $A^\perp \in \rea^{(n-m) \times n}$ a full-rank left-annihilator of $A$. It is rather surprising that the authors of  \cite{Ryalat2018} made such a mistake given that the need of the ``$A^\perp$ condition'' is explicitly stated in \cite[\hspace{-1mm}(6)]{Ryalat2018} for the basic IDA-PBC.
\endrem

\begrem
\lab{rem4} 
Sufficient conditions for \eqref{matching1} to hold are (i) $M_d$ is a constant matrix, and (ii) $G^\perp J_2 M_d^{-1}G=0$. These are the assumptions made in \cite{Donaire2017}, which characterize the class of underactuated mechanical system for which the claim {\bf C1} holds.
\endrem
\subsection{iISS property in Section V.A of \cite{Ryalat2018}}
In  Section V.A of \cite{Ryalat2018} it is claimed the following.\\ 

\noindent {\bf C3.}  The system \eqref{sysdist}, with $d_1=0$, in closed-loop with the control law \cite[\hspace{-1mm}(19)]{Ryalat2018} is iISS with respect to the disturbance $d_2$. 
\\ \ \\
This claim is, unfortunately, incorrect for the following reasons.\\

\noindent {\bf R3.}  To establish the claim an upper bound on the time derivative of the desired energy function $\tilde H(x_q,x_p,x_v)$ given in \cite[\hspace{-1mm}(16)]{Ryalat2018} is computed in \cite[\hspace{-1mm}(25)]{Ryalat2018} as follows
\begin{eqnarray}
\dot{\tilde{H}} &\leq& - \|G^\top M_d^{-1} x_p \|^2_{K_v} + \left(M_d^{-1} x_p \right)^\top d_2. \label{young}
\end{eqnarray}
Then, the authors claim that the bound
\begalis{
& - \|G^\top M_d^{-1} x_p \|^2_{K_v} + \left(M_d^{-1} x_p \right)^\top d_2 \leq \\
&-\frac 12 \lambda_{min}(K_v) \left|G^\top M_d^{-1} x_p \right|^2 +\frac{1}{2 \lambda_{min}(K_v)} \left|d_2\right|^2 
}
follows from the application of Young's inequality, which the authors recall in Section II of \cite{Ryalat2018}. This is, unfortunately, incorrect because---if $G$ is not square---the right hand side of \eqref{young} is not in the form $c_1|y|^2+c_2|y|\,|z|$, with constants $c_1$ and $c_2$. 

The issue here is that the disturbance $d_2$ is {\em not matched}, since it is not in the image of $G$. Indeed, considering a true matched disturbances $d_2=G \hat{d}_2$, \eqref{young} becomes
\begin{align}
\dot{\tilde{H}} \leq& - \left|\left|G^\top M_d^{-1} x_p \right|\right|^2_{K_v} + \left(M_d^{-1} x_p \right)^\top G \, \hat{d}_2 \nonumber \\
=& - \left|\left|G^\top M_d^{-1} x_p \right|\right|^2_{K_v} + \left(G^\top M_d^{-1} x_p \right)^\top  \hat{d}_2 \label{newyoung}
\end{align}
and Young's inequality can be applied.

A further mistake is made in the second equation of \cite[\hspace{-1mm}(26)]{Ryalat2018}, where it is claimed that there exists a $\mathcal{K}_\infty$ function $\alpha$ such that
\begalis{
& -\frac 12 \lambda_{min}(K_v) \left|G^\top M_d^{-1} x_p \right|^2 +\frac{1}{2 \lambda_{min}(K_v)} \left|d_2\right|^2  \\
&\leq -\alpha(|x_p|)+ \frac{1}{2 \lambda_{min}(K_v)} \left|d_2\right|^2.
}
Since $G$ is not square, $\alpha$ cannot be an $\mathcal{K}_\infty$ function of $|x_p|$. Therefore, the iISS property with respect to the input $d_2$ claimed in Section V.A of \cite{Ryalat2018} does not hold.\\

\begrem
\lab{rem5} 
Notice that applying Young's inequality to the correct bound \eqref{newyoung}---corresponding to a {\em bona fide} matched disturbance $\hat d_2$---it follows that 
\begin{align}
\dot{\tilde{H}} \leq& -\frac 12 \lambda_{min}(K_v) \left|y_p \right|^2 +\frac{1}{2 \lambda_{min}(K_v)} \left| \hat{d}_2\right|^2  \label{newalpha}
\end{align}
with $y_p:=G^\top M_d^{-1} x_p$. Then, under an assumption of detectability of $y_p$ when $\hat{d}_2\equiv 0$, the closed-loop is iISS with respect to the disturbance $\hat{d}_2$. It is not surprising that pH systems with damping injection enjoy iISS properties with respect to matched disturbances. 
\endrem

\subsection{Proposition 5.1 of \cite{Ryalat2018}} \label{prop51}
The claims in this proposition are that the system \eqref{sysdist} in closed loop with the control law \cite[\hspace{-1mm}(27)]{Ryalat2018} verifies the following.\\

\noindent {\bf C4.}  The system can be written in the pH form \cite[\hspace{-1mm}(31)]{Ryalat2018}, which is described in new coordinates defined in \cite[\hspace{-1mm}(29)]{Ryalat2018}.

\noindent {\bf C5.} The system is ISS respect to the disturbances $d_1$ and $d_2$.
\\ \ \\
Both claims are, unfortunately, incorrect for the following reasons.\\

\noindent {\bf R4.} Similarly to the mistake made in the claim {\bf C1} above, the ``$G^\perp$ condition''
\begin{eqnarray}
0 \hspace{-2mm}&=&\hspace{-2mm} G^\perp \left\{ M_dM^{-1} \nabla_q H_d - 2 (M^{-1}\mathcal{K} + M_d M^{-1}) \nabla_{x_q} \tilde H   \right. \nonumber \\
&&   - (J_2-R_d) M_d^{-1}(q)p + (J_2-R_d)M_d^{-1}(x_q)p \nonumber \\
&&  \left. - M^{-1} p   + 2 (J_2-R_d)  M_d^{-1} \mathcal{K}  x_v  - 2 M_dM^{-1} x_v \right\}, \nonumber
\end{eqnarray}
needs to be imposed to ensure the matching equation
\begin{eqnarray}
G\,v \hspace{-2mm}&=&\hspace{-2mm} M_dM^{-1} \nabla_q H_d - 2 (M^{-1}\mathcal{K} + M_d M^{-1}) \nabla_{x_q} \tilde H    \nonumber \\
&&  - (J_2-R_d) M_d^{-1}(q)p + (J_2-R_d)M_d^{-1}(x_q)p    \nonumber \\
&& - M^{-1} p + 2 (J_2-R_d)  M_d^{-1} \mathcal{K}  x_v - 2 M_dM^{-1} x_v  \nonumber \\ \label{133}
\end{eqnarray}
admits the solution \cite[\hspace{-1mm}(27)]{Ryalat2018}.  Without this condition, the matching claim is incorrect.

Moreover, the dynamics of $x_q$ obtained from \cite[\hspace{-1mm}(29)]{Ryalat2018} is
\begin{eqnarray}
\dot x_q \hspace{-2mm}&=&\hspace{-2mm} \dot q - \dot x_v \nonumber \\
\hspace{-2mm}&=&\hspace{-2mm} M^{-1}p+d_1 -  M^{-1} \mathcal{K}  \nabla_{x_q}\tilde{H} - \frac 12 M^{-1} p   \nonumber \\   
\hspace{-2mm}&=&\hspace{-2mm}  -  M^{-1} \mathcal{K}  \nabla_{x_q}\tilde{H} + \frac 12 M^{-1} p  +d_1  \nonumber \\   
\hspace{-2mm}&=&\hspace{-2mm}  -  M^{-1} \mathcal{K}  \nabla_{x_q}\tilde{H} + M^{-1} x_p  - M^{-1} \mathcal{K} x_v + d_1 \nonumber \\   
\hspace{-2mm}&=&\hspace{-2mm}  -  M^{-1} \mathcal{K}  \nabla_{x_q}\tilde{H} + M^{-1}M_dM_d^{-1} x_p  - M^{-1} \mathcal{K} x_v + d_1 \nonumber \\   
\hspace{-2mm}&=&\hspace{-2mm}  -  M^{-1} \mathcal{K}  \nabla_{x_q}\tilde{H} + M^{-1}M_d \nabla_{x_p}\tilde{H}  - M^{-1} \mathcal{K} \nabla_{x_v}\tilde{H} + d_1 \nonumber    
\end{eqnarray}
which shows that the first row of the desired dynamics in \cite[\hspace{-1mm}(31)]{Ryalat2018} is not achieved.\\

\noindent {\bf R5.}  The claim of ISS is incorrect because, on one hand, the closed-loop dynamics is not in the form \cite[\hspace{-1mm}(31)]{Ryalat2018}---for the reason given above. On the other hand, as in the rebuttal of {\bf C3}, it is easy to see that Young's inequality is erroneously used to get the first bound in \cite[\hspace{-1mm}(32)]{Ryalat2018}. Moreover, the claim that the function $\alpha$, appearing in the second bound, is $\mathcal K_\infty$ for $|x_q,x_p,x_v|$ is wrong because $G$ is not square.

\section{Examples of \cite{Ryalat2018}}
In this section we prove that the two examples considered in \cite{Ryalat2018} are incorrect.
\subsection{The inertia wheel pendulum}
The first example proposed in \cite{Ryalat2018} is the inertia wheel pendulum (IWP). The controller for this example is designed using  the result in Proposition 5.1. We will show next that, as shown in Subsection \ref{prop51} for the general case, the proposed controller for the IWP does not satisfy the matching equation.

Consider the following matrices and functions for the IWP \cite{Ryalat2018}:
\begin{eqnarray}
M&=&\left[ \begin{array}{cc} k_1 & k_2 \\ k_2 & k_2 \end{array} \right], \hspace{1.5mm} G=\left[ \begin{array}{c} 0 \\ 1 \end{array} \right], \quad V=k_3[1+\cos(q_1)],  \nonumber \\
  M_d&=&\Delta\left[ \begin{array}{cc} m_1 & m_2 \\ m_2 & m_3 \end{array} \right]\nonumber \\
   V_d&=&-k_3\gamma_1 \cos(q_1) + \frac 12 K_p [\epsilon k_1 \gamma_1 q_1 +q_2]^2, \nonumber \\
&& \hspace{-15mm} \tilde{V}=-k_3\gamma_1 \cos(x_{q_1}) + \frac 12 K_p [\epsilon k_1 \gamma_1 x_{q_1} +x_{q_2}]^2. \nonumber
\end{eqnarray}
From \eqref{133}, the control law should satisfy
\begin{eqnarray}
G\,v &=& M_dM^{-1} \nabla_q V_d - 2 M^{-1}\mathcal{K} \nabla_{x_q} \tilde V - 2 M_d M^{-1} \nabla_{x_q} \tilde V  \nonumber \\
&&  - M^{-1} p - 2 \left( R_d M_d^{-1} \mathcal{K}  x_v +  M_dM^{-1} \right) x_v, \label{matchexamp}
\end{eqnarray}
with $\mathcal{K} = GK_iG^\top$, $R_d=GK_vG^\top$ and $x_q=q- x_v$.
\\ \ \\
The terms in the right-hand-side of \eqref{matchexamp} can be computed as follows
\begin{eqnarray}
 M_dM^{-1} \nabla_q V_d&= & \left[ \begin{array}{cc} 1  \\ 1 \end{array} \right] K_p \epsilon \mathcal{S} (\epsilon \gamma_1 k_1 q_1 + q_2) \nonumber \\
 &+& \gamma_1 k_2 k_3 (m_1-m_2) \sin(q_1) \label{term1} \\
-2 M^{-1}\mathcal{K} \nabla_{x_q} \tilde V  &\hspace{-6mm}  =& \hspace{-6mm} \left[ \begin{array}{c} 2k_2  \\ -2k_1 \end{array} \right] \frac{K_i K_p}{(k_1-k_2)k_2} (\epsilon \gamma_1 k_1x_{q_1} + x_{q_2}), \nonumber \\
\label{term2} \\
 - 2 M_d M^{-1} \nabla_{x_q} \tilde V &=& - \left[ \begin{array}{c} 2  \\ 2 \end{array} \right] K_p \epsilon \mathcal{S}[\epsilon \gamma_1 k_1  x_{q_1} +  x_{q_2} ]  \nonumber \\ 
&&  + \gamma_1 k_2 k_3 (m_1-m_2) \sin(x_{q_1} ), \label{term3} \\
 -M^{-1} p &=& \frac{1}{k_2(k_1-k_2)}  \left[ \begin{array}{cc} -k_2 (p_1-p_2)  \\ k_2 p_1 +k_1 p_2 \end{array} \right], \label{term4} \\
&&  \hspace{-3.4cm}-2 \left( R_d M_d^{-1} \mathcal{K}  x_v +  M_dM^{-1} \right) x_v = - 2 \mathcal{L} \nonumber \\
 && - \hspace{-.5cm} \left[ \begin{array}{c} 0  \\  1 \end{array} \right]  \frac{2 K_i m_1 K_v x_{v_2}}{k_1(k_1-k_2)(m_1m_3-m_2)^2}. \nonumber \\ \label{term5}
\end{eqnarray}
with 
\begin{align}
\mathcal{S}:=&1+ \gamma_1 k_1 k_2 (m_1-m_2)\nonumber \\
\mathcal{L}:=&\left[ \begin{array}{c} k_2 (m_1-m_2)x_{v_1} + k_1 (m_2-m_1) x_{v_2}  \\   k_2 (m_1-m_3)x_{v_1} + k_1 (m_3-m_1) x_{v_2}   \end{array} \right].  \nonumber
\end{align}
Using \eqref{term1}-\eqref{term5} in \eqref{matchexamp}, we obtain that the matching equation \eqref{matchexamp} has the form
\begin{eqnarray}
\left[ \begin{array}{c} 0  \\  1 \end{array} \right]  v &=& \left[ \begin{array}{c} (\star)  \\  (\star \star) \end{array} \right],   \label{nomatch}
\end{eqnarray}
where $(\star)$ and $(\star\star)$ are non-zero. From \eqref{nomatch} we conclude that the matching equation cannot be satisfied.
\subsection{Rotary inverted pendulum}
A controller for the rotary inverted pendulum (RIP) is designed using Proposition 4.1 in \cite{Ryalat2018}. We will shown in this section that the proposed controller does not satisfy the matching equation and therefore the stability claims have no theoretical support.

The control law for the RIP is obtained from \cite[\hspace{-1mm}(15)]{Ryalat2018} using the following matrices:
 \begin{eqnarray}
&&   M_d^{-1}= \left[ \begin{array}{cc} \frac{\Delta m_3}{\Delta_d} & \frac{\Delta \sigma \cos(q_1) (\cos(q_1) +\epsilon)}{\Delta_d \gamma} \\ \frac{\Delta \sigma \cos(q_1) (\cos(q_1) +\epsilon)}{\Delta_d \gamma}  & \frac{\Delta (\cos(q_1) +\epsilon)}{\Delta_d } \end{array} \right], \nonumber \\
&&  \nabla_{q_i} M_d^{-1} = \left[ \begin{array}{cc} \mathcal{B}_1 & \mathcal{B}_2 \\  \mathcal{B}_2  & \mathcal{B}_3 \end{array} \right], \nonumber
\end{eqnarray}
where all the functions and parameters are defined in Section VII of \cite{Ryalat2018}.

From \cite[\hspace{-1mm}(22)]{Ryalat2018}, the controller should satisfy
\begin{align}
G\,v =& \frac 12 M_dM^{-1} \sum e_i p^\top \nabla_{q_i} M_d^{-1} p + (J_2-R_d) M_d^{-1} \mathcal{K}  x_v  \nonumber \\
& - \frac 12 M_dM^{-1} \sum e_i x_p^\top \nabla_{q_i} M_d^{-1} x_p - \mathcal{K}  \dot{x}_v.  \label{matchexamp2}
\end{align}
with $x_p=p+\mathcal{K} x_v$.

The terms on right-hand-side of \eqref{matchexamp2} can be computed as follows
{\small
\begin{eqnarray}
&& \hspace{-8mm} \frac 12 M_dM^{-1} \sum e_i p^\top \nabla_{q_i} M_d^{-1} p = \left[ \begin{array}{cc} \frac 12 \mathcal{B}_1 p_1^2 + \mathcal{B}_2 p_1 p_2 +\frac 12  \mathcal{B}_2 p_2^2  \\ 0 \end{array} \right],  \nonumber \\\label{term12} \\
&&  (J_2-R_d) M_d^{-1} \mathcal{K}  x_v  =  {\mathcal Q} x_{v_2},\nonumber \\ \label{term22} \\
&& \frac 12 M_dM^{-1} \sum e_i x_p^\top \nabla_{q_i} M_d^{-1} x_p  = \left[ \begin{array}{cc} {\mathcal P}  \\ 0 \end{array} \right],  \label{term32} \\
&&  - \mathcal{K}  \dot{x}_v =  \left[ \begin{array}{c} 0  \\  \dot{x}_{v_2} \end{array} \right], \label{term42} 
\end{eqnarray}
}
where we defined
\begalis{
{\mathcal Q} & :=\frac 12 \mathcal{B}_1 x_{p_1}^2 + \mathcal{B}_2 x_{p_1} x_{p_2} +\frac 12  \mathcal{B}_2 x_{p_2}^2 \\
{\mathcal P} &:= \begmat{  j_2 \frac{\Delta (\cos(q_1) +\epsilon)}{\Delta_d }    \\  -j_2 \frac{\Delta \sigma \cos(q_1) (\cos(q_1) +\epsilon)}{\Delta_d \gamma} - k_v \frac{\Delta (\cos(q_1) +\epsilon)}{\Delta_d }}. 
}

Using \eqref{term12}-\eqref{term42} in \eqref{matchexamp2}, we show that the matching equation \eqref{matchexamp2} has the form
\begin{eqnarray}
\left[ \begin{array}{c} 0  \\  1 \end{array} \right]  v &=& \left[ \begin{array}{c} (\star)  \\  (\star \star) \end{array} \right],   \label{nomatch2}
\end{eqnarray}
where $(\star)$ and $(\star\star)$ are non-zero. Indeed,
\begin{align}
(\star) =& \frac 12 \mathcal{B}_1 p_1^2 + \mathcal{B}_2 p_1 p_2 +\frac 12  \mathcal{B}_2 p_2^2 + j_2 \frac{\Delta (\cos(q_1) +\epsilon)}{\Delta_d } x_{v_2}  \nonumber \\
& - \mathcal{B}_2 p_1 (p_2+x_{v_2}) - \frac 12  \mathcal{B}_2 (p_2+x_{v_2})^2 - \frac 12 \mathcal{B}_1 p_1^2 \nonumber \\
=&   j_2 \frac{\Delta (\cos(q_1) +\epsilon)}{\Delta_d } x_{v_2}  - \mathcal{B}_2 p_1 x_{v_2} - \mathcal{B}_2 p_2x_{v_2}  - \frac 12  \mathcal{B}_2 x_{v_2}^2 \nonumber
\end{align}
which clearly shows, together with \eqref{nomatch2}, that the matching equation cannot be satisfied.

Also, notice that the terms \eqref{term12}  and \eqref{term32} are conspicuous by their absence in the controller computed in \cite[\hspace{-1mm}(48)]{Ryalat2018}.

\section{Conclusions}

We have shown that the claims made in Proposition 4.1, Section V.A of \cite{Ryalat2018} and Proposition 5.1, as well as both examples presented there, are wrong. There are many other examples in the literature, {\em e.g.}, \cite{KWAKRE,KWAONU,ZENetal} where PBC has been erroneously applied. This situation is not critical in application-oriented publications, where the main emphasis is not in mathematical correctness but on the fact that the proposed design yields a satisfactory performance, validated by experiments. It is however critical in the present case, where the publication is made in a theoretical journal, where mathematical rigor is of prime importance. Moreover, it is very harmful to the control community, to allow the publication of  a paper that claims to extend a well-established method. Particularly considering that there are already several schemes that achieve this objective, which are cited in \cite{Ryalat2018}, that is  \cite{Donaire2009,Romero2013,Donaire2017}. 



\begin{thebibliography}{9}

\bibitem{Ryalat2018}
M.~Ryalat and D.~S.~Laila, ``A robust IDA-PBC approach for handling uncertainties in underactuated mechanical systems,''  \emph{IEEE Transactions on Automatic Control}, vol. 63, no. 10, pp. 3495-3502, 2018.


\bibitem{Donaire2017}
A.~Donaire, J.~G.~Romero, R.~Ortega, B.~Siciliano and M.~Crespo, ``Robust IDA-PBC for underactuated mechanical systems subject to matched disturbances,'' \emph{International Journal of Robust and
Nonlinear Control}, vol.~27, no.~6, pp. 1000-1016, 2017.
 
\bibitem{KWAKRE}
A. Kwasinski and P. Krein, ``Passivity-Based Control of Buck Converters with Constant-Power Loads,'' in {\em Proc. Power Electron. Spec. Conf.}, pp.259-265, Jun. 2007.

\bibitem{KWAONU}
A. Kwasinski and C. N. Onwuchekwa, ``Dynamic behavior and stabilization of DC microgrids with instantaneous constant-power loads,'' {\em IEEE Trans. on Power Electron.}, vol.26, no.3, pp.822-834, Mar., 2011.

\bibitem{ORTtac}
R.~Ortega, M.~Spong, F.~Gomez-Estern, and G.~Blankenstein, ``Stabilization of a
  class of underactuated mechanical systems via interconnection and damping
  assignment,'' \emph{IEEE Transactions on Automatic Control}, vol.~47, no.~8,
  pp. 1218-1233, 2002.


\bibitem{ZENetal}
J. Zeng, Z. Zhang and W. Qiao, ``Interconnection and Damping Assignment Passivity-Based Controller for a DC- DC Boost Converter With a Constant Power Load'', {\em IEEE Trans. Ind. Appl.}, vol. 50, no. 4, pp.2314-2322, Jul., 2014. 	

\bibitem{Donaire2009}
A.~Donaire and S. Junco, ``On the addition of integral action to port- controlled Hamiltonian systems,'' \emph{Automatica}, vol.~45, no.~8, pp. 1910-1916, 2009.

\bibitem{Romero2013}
J.~G.~Romero, A.~Donaire and R.~Ortega, ``Robust energy shaping control of mechanical systems,'' \emph{Systems \& Control Letters}, vol.~62, pp. 770-780, 2013.


\end{thebibliography}
\end{document}